\documentclass[letterpaper,10pt]{article}
\usepackage[textwidth=11.5cm, textheight=22cm, top=2.0cm, centering]{geometry}
\usepackage{setspace}
\usepackage[T1]{fontenc}
\usepackage{lmodern}
\usepackage[english]{babel}
\usepackage{microtype} 
\usepackage[scaled=0.95]{inconsolata} 
\makeatletter
\renewcommand{\texttt}[1]{{\ttfamily #1}}

\renewcommand{\mathtt}[1]{\text{\texttt{#1}}} 
\makeatother
\usepackage{csquotes}
\usepackage{etoolbox}
\usepackage{bbding}
\usepackage{changepage}
\usepackage{booktabs}
\usepackage{graphicx}
\usepackage[x11names, dvipsnames]{xcolor}
\definecolor{Linkz}{RGB}{30, 110, 170}
\definecolor{Darkenta}{RGB}{185, 35, 90}
\definecolor{Lightenta}{RGB}{254, 232, 255}
\definecolor{Reference}{RGB}{35, 180, 90}
\definecolor{Periwinkle}{RGB}{102, 51, 255}
\definecolor{Greeno}{RGB}{0, 140, 100}
\definecolor{Leeno}{RGB}{239, 255, 232}
\usepackage{amsmath, amsthm, amscd, amssymb, stmaryrd}
\usepackage{stackengine}
\usepackage[square]{natbib}

\usepackage[
colorlinks=true,
linkcolor=Darkenta,
citecolor=Greeno,
urlcolor=Darkenta,
]{hyperref}
\usepackage{titlesec} 
\usepackage{abstract} 
\usepackage{enumitem} 
\usepackage{ccicons} 
\usepackage{float}
\usepackage{placeins} 
\usepackage{tikz} 
\usepackage{tikz-cd}
\usetikzlibrary{positioning,fit}
\usetikzlibrary{arrows.meta} 
\usepackage{adjustbox} 
\usepackage{pgfmath} 
\usepackage{ifthen} 

\newtheoremstyle{upright}
{6pt plus 2pt minus 2pt} 
{6pt plus 2pt minus 2pt} 
{\normalfont} 
{} 
{\bfseries} 
{.} 
{.5em} 
{} 
\theoremstyle{upright}
\theoremstyle{upright}
\newtheorem{theorem}{Theorem}[section]
\newtheorem{remark}[theorem]{Remark}

\newtheorem{definition}[theorem]{Definition}

\newtheorem{proposition}[theorem]{Proposition}
\newtheorem{lemma}[theorem]{Lemma}

\newtheorem{corollary}[theorem]{Corollary}

\makeatletter
\renewenvironment{proof}[1][Proof]{%
	\par\pushQED{\qed}%
	\normalfont
	\topsep6\p@\@plus6\p@\relax
	\trivlist
	\item[\hskip\labelsep\slshape #1\@addpunct{.}]%
}{%
	\popQED\endtrivlist\@endpefalse
}

\makeatother

\usepackage{algpseudocode}
\usepackage[most]{tcolorbox}
\usepackage{listings}
\newtcolorbox{breakbox}[2][]{%
	breakable,
	={#2},
	fonttitle=\bfseries,
	colback=white,
	colframe=black!20,
	coltitle=black,
	colbacktitle=white,
	boxrule=0.5pt,
	arc=0pt,
	boxsep=7pt,
	left=3pt,
	right=2pt,
	top=2pt,
	bottom=4pt,
	fontupper=\small\sffamily, 
	#1
}

\AtBeginEnvironment{quotation}{\sffamily}
\renewenvironment{quotation}
{\small\vspace{0.5em}\begin{adjustwidth}{4em}{4em}%

		\setlength{\parindent}{0pt}%
		\setlength{\parskip}{\medskipamount}%
	}
	{\end{adjustwidth}\vspace{1em}}
\newcommand{\customsectionstyle}[2]{%
	\titleformat{\section}[block]
	{\normalfont\fontsize{#1}{1.2\dimexpr#1\relax}\selectfont\centering}
	{\thesection}{1em}%
	{%
		\ifthenelse{\equal{#2}{true}}{\MakeUppercase}{\relax}%
	}%
}
\newcommand{\customsectionspacing}[3]{%
	\titlespacing*{\section}{#1}{#2}{#3}%
}
\newcommand{\customsubsectionstyle}[2]{%
	\titleformat{\subsection}[block]
	{\normalfont\fontsize{#1}{1.2\dimexpr#1\relax}\selectfont\centering}
	{\thesubsection}{1em}%
	{%
		\ifthenelse{\equal{#2}{true}}{\MakeUppercase}{\relax}%
	}%
}
\newcommand{\customsubsectionspacing}[3]{%
	\titlespacing*{\subsection}{#1}{#2}{#3}%
}
\customsectionstyle{14pt}{true}
\customsectionspacing{0pt}{3.5ex plus 1ex minus 0.2ex}{2.5ex}
\customsubsectionstyle{11pt}{true}
\customsubsectionspacing{0pt}{4ex plus 1ex minus 0.2ex}{2ex}
\usepackage{caption}

\makeatletter
\newcommand{\shorttitle}[1]{\def\@shorttitle{#1}}
\newcommand{\email}[1]{\def\@email{#1}}
\newcommand{\metadata}[1]{\def\@metadata{#1}}
\renewcommand{\maketitle}{%
	\begin{center}
		\vfill
		{\fontsize{18pt}{19pt}\selectfont \@title \par}
		\vspace{1em}
		{\normalsize \@author \par}
		\vspace{0.1em}
		{\normalsize \@date \par}
	\end{center}
}
\makeatother

\begin{document}
\newcommand{\Z}{\mathbb{Z}}
\newcommand{\N}{\mathbb{N}}
\newcommand{\Q}{\mathbb{Q}}

\newcommand{\mult}{\mathtt{Mult}}

\newcommand{\HtenZ}{{\rm H}\Z}
\newcommand{\HtenN}{{\rm H}\N}
\newcommand{\HtenQ}{{\rm H}\Q}
\newcommand{\Hten}[3]{\mathrm{H}{#1}(#2,#3)}
\newcommand{\D}{\mathsf{D}}
\newcommand{\U}{\mathsf{U}}
\newcommand{\TheoryT}{\mathrm{T}}

\newcommand{\Prov}{\operatorname{Prov}}
\newcommand{\ProvRA}{\operatorname{Prov}_{\TheoryT}}
\newcommand{\ProofRA}{\operatorname{Proof}_{\TheoryT}}
\newcommand{\FormPred}{\operatorname{Form}}
\newcommand{\Sent}{\operatorname{Sent}}
\newcommand{\SigmaOneSent}{\operatorname{Sigma1Sent}}
\newcommand{\SatSigmaOne}{\operatorname{Sat}_{\Sigma_1}}
\newcommand{\DiophThree}{\operatorname{Dioph}_3}
\newcommand{\RFNSigmaOne}{\operatorname{RFN}_{\Sigma_1}}


\title{\uppercase{Considering The\\
Satisfiability of Cubic\\
Diophantine Equations}}
\author{Milan Rosko}
\date{April 2026}

\maketitle

\begin{abstract}
\vspace{-1.5ex}\footnotesize{
	Our contribution is a bounded cubic compilation theorem. For each fixed resource parameter $k$, syntactic proof checking at resource level $k$ is faithfully represented by a finite bounded-domain system of cubic polynomial equations. Every emitted equation has degree at most 3. Degree-3 terms arise only when a linear selector variable activates a quadratic verification obligation.

	Earlier versions of this manuscript claimed a reduction from unbounded theoremhood to satisfiability of a fixed bounded-domain cubic polynomial instance. That claim is withdrawn. The error and its source are identified precisely. The bounded construction, the degree bookkeeping, and the \textsc{Zeckendorf}-based carryless encoding stand independently of the withdrawn claim.

	The note closes by identifying the uniformization gap that separates a family of decidable bounded slices from a single \emph{many-one reduction} target, and records why closing that gap would require a compression principle not supplied here.
	}
\end{abstract}

\section{Introduction}

\begin{remark}

	\textsc{Hilbert's Tenth Problem} \citep{hilbert1900} is, in a sense, the question of how far one can jump before jumping becomes impossible. \textsc{MRDP} \citep{rdp61,matiyasevich70,davis1973,mrdpbook} settled that the other side of the street is unreachable. Conversely, the near ground is solid. But between the two lies a stretch where the takeoff itself starts to wobble, as the question of whether we could jump becomes entangled with the jump itself. The present work was the product of a year-long occupation with this question, until we fell backward.

	One could also say that we missed the high bar of \emph{proof theory}: namely, to \emph{fail elegantly}.

\end{remark}
\begin{definition}[\textsc{Gasarch} notation]
\label{def:hn-notation}
	We fix a convenient \citep{gasarch21} notation for parameterized Diophantine solvability throughout. For $S\in\{\N,\Z\}$ and $d,n\in\N$, let
	\begin{equation}
		\Hten{S}{d}{n}
	\end{equation}
	denote the decision problem whose inputs are integer-coefficient polynomials $P$ of total degree at most $d$ in at most $n$ variables, and whose question is whether
	\begin{equation}
		\exists \vec{x}\in S^m\; P(\vec{x})=0
		\qquad(m\leq n).
	\end{equation}
	We write $\Hten{S}{d}{n}=\D$ when this problem is decidable and $\Hten{S}{d}{n}=\U$ when it is undecidable.
\end{definition}

\begin{lemma}[Systems and aggregation]
\label{prop:trivial-systems}
	A finite system
	\begin{equation}
		P_1(\vec{x})=0,\ldots,P_m(\vec{x})=0
	\end{equation}
	over $\N$ or $\Z$ can always be reduced to one equation by
	\begin{equation}
		\sum_{i=1}^m P_i(\vec{x})^2=0.
	\end{equation}
	This doubles the degree. Hence ordinary sum-of-squares aggregation is unsuitable when one wants to preserve a cubic degree bound.
\end{lemma}

\begin{remark}
	The original motivation was to compile proof checking into polynomial constraints with local degree control. Syntactic verification, formula decoding, target equality, and inference checking are essentially quadratic. Cubic terms appear only when a linear selector activates a quadratic obligation. The corrected result retains this bounded degree analysis in full. What it does not do is treat the resulting bounded family as a fixed unbounded universal instance.
\end{remark}

\begin{definition}[Object theory]
\label{def:register-arithmetic}

	The theoremhood predicate is taken for a fixed recursively axiomatized undecidable theory $\TheoryT$ in the language of arithmetic
		\begin{equation}
			\mathcal{L} = \{0,S,+,\times,=\}.
		\end{equation}
	For concreteness one may take $\TheoryT$ to be Robinson arithmetic $Q$, or any consistent recursively axiomatized extension of $Q$. The undecidability of theoremhood for such theories is classical \citep{kleene52,boolos07}.

	The bounded verifier below uses only the recursive presentation of the axioms and rules of $\TheoryT$.
\end{definition}

\begin{remark}
	In elementary examples the arithmetic part contains the successor and addition axioms
	\begin{align}
		&S(x) \ne 0, &
		&S(x)=S(y)\to x=y, &
		&x\ne 0\to \exists y\;x=S(y),\\
		&x+0=x, &
		&x+S(y)=S(x+y),
	\end{align}
	together with bounded induction for $\Delta_0$-formulae and $\mathrm{B}\Sigma_1$-collection:
	\begin{equation}
		(\mathrm{I}\Delta_0+\mathrm{B}\Sigma_1)\upharpoonright\mathcal{L}.
	\end{equation}
\end{remark}

\begin{remark}
	The checker only inspects finite formula and proof codes. Any auxiliary arithmetic trace used by the encoding is verified by bounded inspection; the undecidability used in the corrigendum comes from theoremhood of the chosen theory $\TheoryT$.
\end{remark}

\begin{remark}
	The degree threshold arises from the difference between checking a supplied finite witness and asserting the existence of an unbounded witness. Degree-$2$ constraints handle equality, booleanity, structural decoding, and many local checks. Degree-$3$ appears when a selector variable, of degree $1$, activates a quadratic verification condition.
\end{remark}

\begin{lemma}[Guard correctness and degree]
\label{lem:guard-correctness}
	Let $E(\vec{x})$ be a polynomial of degree $\leq 2$ and let $s$ be a selector variable. Then
	\begin{equation}
		s\cdot E(\vec{x})=0
	\end{equation}
	has degree at most $3$. In the emitted systems, selectors are used only to activate quadratic obligations; no selector is multiplied by a cubic expression.
\end{lemma}

\begin{remark}
	We recall the degree hierarchy of Diophantine equations over $\mathbb{N}$:
	\begin{enumerate}[label=\textnormal{(\roman*)}]
		\item \emph{Degree 2}: Certain structured quadratic forms are decidable via classical number-theoretic methods (\textsc{Hasse--Minkowski} for quadratic forms over $\mathbb{Q}$, \textsc{Lagrange} \citep{lagrange1770} for sums of squares). Degree-$2$ constraints lack sufficient expressive power to encode arbitrary proof checking.
		\item \emph{Degree 4}: \citet{jones80,jones82} showed that every recursively enumerable set is the solution set of a degree-$4$ Diophantine equation with $58$ variables. Undecidability follows classically from \textsc{MRDP} \citep{rdp61,matiyasevich70,davis1973,mrdpbook}.
		\item \emph{Degree 3}: Black-box reductions from \textsc{MRDP} composed with \textsc{Jones Shielding} \citep{jones80} and a quartic-to-cubic-system closure yield $\Sigma_1$-completeness of cubic \emph{systems} (with unbounded variables), but furnish no degree-local construction and no single-equation result at degree~$3$.
	\end{enumerate}
	\citet{gasarch21} summarizes the state of the problem as one in which progress on $\Hten{\N}{d}{n}=\U$ or $\Hten{\Z}{d}{n}=\U$ appears to have stalled \citep{gasarch21}. \citet{jones80} explicitly singled out degree~$3$ over $\mathbb{N}$ as the sole open case between the understood quadratic setting and high-degree universality constructions \citep{jones80}. The corrected manuscript contributes to the bounded side of that problem without resolving the single-equation question.
\end{remark}

\section{Corrigendum}
\label{sec:error}

\begin{remark}
	Let $\ProvRA(t)$ denote unbounded theoremhood, or equivalently the existence of an unrestricted proof of the formula coded by $t$. The earlier versions asserted a \emph{many-one reduction}
	\begin{equation}
		\ProvRA(t)
		\quad\Longleftrightarrow\quad
		\mathrm{Sat}_{B_{f(t)}}(t)
	\end{equation}
	for some computable function $f$, where $\mathrm{Sat}_{B_k}(t)$ denotes satisfiability of the bounded-domain polynomial slice at resource level $k$. This claim is incorrect.

	The inference fails because it requires a computable, uniform selection of a sufficient resource bound from the input $t$ alone. Each fixed bounded slice $\mathrm{Sat}_{B_k}$ is decidable by finite search. If such a selector $f$ existed, one could decide the unbounded theoremhood predicate by computing $f(t)$ and running finite search over the corresponding slice. Since $\ProvRA$ is recursively enumerable but not decidable for the chosen theory $\TheoryT$, no such selector can exist. The construction supplies no mechanism for producing $f$, so the passage from the bounded family to the reduction is unjustified.
\end{remark}

\begin{remark}
	The correct statement is slice-wise. For each fixed resource parameter $k$,
	\begin{equation}
		\operatorname{BoundedProof}_k(t)
		\quad\Longleftrightarrow\quad
		\mathrm{Sat}_{B_k}(t).
	\end{equation}
	Both predicates include the size convention $t\leq B_k$ and are false when that convention fails.
	Unbounded theoremhood is recovered only as the existential union over all slices,
	\begin{equation}
		\mathrm{Sat}^{\cup}(t) :\iff \exists k\;\mathrm{Sat}_{B_k}(t),
	\end{equation}
	which is a valid recursively enumerable predicate but not the satisfiability predicate of any single fixed bounded-domain polynomial instance.
\end{remark}

\begin{remark}
	A \emph{many-one reduction} is a single computable translation
	\begin{equation}
		t\longmapsto x(t),
		\qquad
		A(t)\iff B(x(t)).
	\end{equation}
	A bounded approximation more often gives
	\begin{equation}
		A(t)\iff \exists k\;B_k(t),
	\end{equation}
	where each $B_k$ is decidable. Such a statement proves recursive enumerability, but it does not supply the single target instance required.

	This is the substantive gap. Bounded correctness and unbounded completeness are separated by a uniformity problem: the resource parameter must either be internalized into one Diophantine predicate or selected by a computable function. The present construction establishes the bounded side only, and this is stated explicitly going forward.
\end{remark}

\section{Carryless Infrastructure}
\label{sec:arith-infra}

\begin{theorem}[\textsc{Zeckendorf} representation]
\label{thm:zeckendorf}
	Every $n>0$ has a unique representation
	\begin{equation}
		n=\sum_{k\in S}F_k,
	\end{equation}
	where $S\subseteq\{2,3,4,\ldots\}$ is finite and contains no consecutive indices \citep{zeckendorf}.
\end{theorem}

\begin{theorem}[\textsc{Carryless Pairing}]
\label{def:carryless-pairing}
	Let $Z(x)$ denote the \textsc{Zeckendorf} support of $x$. Define two separated support bands
	\begin{align}
		E(x)&:=\{2k:k\in Z(x)\},\\
		O(x,y)&:=\{B(x)+(2j-1):j\in Z(y)\},
	\end{align}
	where $B(x)$ is chosen beyond the even support of $x$. The pair code is
	\begin{equation}
		\mathrm{pair}(x,y):=F_2+\sum_{i\in E(x)\cup O(x,y)}F_i.
	\end{equation}
	The construction and its primitive recursive projections are developed in detail in our earlier work \citep{rosko2025fibonacci}.

	The function $\mathrm{pair}:\N^2\to\N$ is injective and has primitive recursive projections $\mathrm{hd}$ and $\mathrm{tl}$ such that
	\begin{equation}
		\mathrm{hd}(\mathrm{pair}(x,y))=x,
		\qquad
		\mathrm{tl}(\mathrm{pair}(x,y))=y.
	\end{equation}
\end{theorem}

\begin{remark}
	The supports of $x$ and $y$ occupy disjoint Fibonacci bands. By \citet{zeckendorf}, the union support can be recovered and split into its two components without any positional carry interaction. Projections are therefore realised by bounded search, which is what keeps the downstream constraints quadratic.
\end{remark}

\begin{definition}[List encoding and line extraction]
\label{def:list-code}
	Define
	\begin{align}
		\mathrm{code}_{\mathrm{list}}([])&:=0,\\
		\mathrm{code}_{\mathrm{list}}(x_0,\ldots,x_{n-1})&:=\mathrm{pair}(x_0,\mathrm{code}_{\mathrm{list}}(x_1,\ldots,x_{n-1})).
	\end{align}
	The $i$-th line is obtained by primitive recursion:
	\begin{align}
		\mathrm{tail}^0(s)&:=s,\\
		\mathrm{tail}^{i+1}(s)&:=\mathrm{tl}(\mathrm{tail}^i(s)),\\
		\mathrm{line}(s,i)&:=\mathrm{hd}(\mathrm{tail}^i(s)).
	\end{align}
\end{definition}

\begin{remark}
	For a fixed bounded slice, line extraction is represented by a finite chain of witnesses. This keeps the local equations quadratic before selector activation. It does not, by itself, give a constant-overhead Diophantine predicate for an unbounded run or proof history---a further reason why the slice-wise result does not extend to the unbounded setting for free.
\end{remark}

\begin{definition}[Proof predicate]
\label{def:proof-predicate}
	A proof code $p$ for a target $t$ is a finite sequence of formula codes such that each line is either an axiom instance or follows by modus ponens from earlier lines, and the final line is $t$. We write
	\begin{equation}
		\ProofRA(p,t)
	\end{equation}
	for this primitive recursive relation, and
	\begin{equation}
		\ProvRA(t):\iff\exists p\;\ProofRA(p,t).
	\end{equation}
	The arithmetization follows standard lines; see \citet{hajekpudlak93,kleene52,feferman60}.
\end{definition}

\begin{definition}[Bounded proof predicate]
\label{def:bounded-proof-predicate}
	For each resource parameter $k$, fix a computable bound $B_k$ for the proof length, formula-code size, proof-code size, and auxiliary witness values admitted in the $k$-slice. Define
	\begin{equation}
		\operatorname{BoundedProof}_k(t)
	\end{equation}
	to mean that $t\leq B_k$ and that there is a $\TheoryT$-proof of the formula coded by $t$ satisfying those resource restrictions. If $t>B_k$, then $\operatorname{BoundedProof}_k(t)$ is false by definition.
\end{definition}

\begin{lemma}[Primitive recursiveness]
\label{lem:syntactic-prim-rec}
	The syntactic predicates for formulahood, sentencehood, proof checking, axiom matching, and modus ponens verification are primitive recursive.
\end{lemma}

\begin{remark}
	All checks are bounded inspections of finite codes. Axiom matching is structural pattern matching, and modus ponens witnesses are searched only among earlier lines. No unbounded quantifier is involved; see \citet{hajekpudlak93,boolos07} for the bounded arithmetic setting.
\end{remark}

\begin{lemma}[Basic constraints]
\label{lem:basic-constraints}
	Over $\N$:
	\begin{align}
		x=y &\iff (x-y)^2=0,\\
		b\in\{0,1\} &\iff b(1-b)=0.
	\end{align}
	Both are degree $2$ constraints.
\end{lemma}

\begin{definition}[Finite bounded graph expansion]
\label{def:bounded-graph-expansion}
	Let $R\subseteq[0,B]^m$ be a finite relation. The bounded graph expansion
	\begin{equation}
		\operatorname{Graph}_{B}(R;z_1,\ldots,z_m)
	\end{equation}
	is the following finite polynomial system. Introduce one selector $e_{\vec{a}}$ for each tuple $\vec{a}=(a_1,\ldots,a_m)\in R$, and impose
	\begin{align}
		e_{\vec{a}}(1-e_{\vec{a}})&=0 &&(\vec{a}\in R),\\
		\sum_{\vec{a}\in R}e_{\vec{a}}-1&=0,\\
		z_j-\sum_{\vec{a}\in R}a_j e_{\vec{a}}&=0 &&(1\leq j\leq m).
	\end{align}
	All emitted equations have degree at most $2$.

	For a bounded primitive-recursive function $F:[0,B]^r\to\N$, the notation
	\begin{equation}
		y=F(x_1,\ldots,x_r)
	\end{equation}
	inside the construction means this graph expansion for the finite relation
	\begin{equation}
		R^F_B=\{(x_1,\ldots,x_r,y)\in[0,B]^{r+1}: y=F(x_1,\ldots,x_r)\}.
	\end{equation}
	This is the bounded graph of $F$ inside the slice; if the true value of $F$ lies above $B$, the corresponding bounded constraint has no satisfying tuple.
	Thus the emitted polynomial system contains only constants, variables, additions, and multiplications; it contains no primitive-recursive function symbol. In particular,
	\begin{equation}
		\operatorname{Pair}_B(z,x,y),\quad
		\operatorname{Hd}_B(z,x),\quad
		\operatorname{Tl}_B(z,y),\quad
		\operatorname{Line}_{B,i}(s,\ell)
	\end{equation}
	are abbreviations for finite systems obtained from the bounded graphs of $\mathrm{pair}$, $\mathrm{hd}$, $\mathrm{tl}$, and the fixed-index line-extraction map, with the graph constants effectively computed from the \textsc{Zeckendorf} definitions before the polynomial system is emitted.

	If such an expansion occurs inside a selected obligation, it is not inserted unguarded. The notation
	\begin{equation}
		s\cdot\operatorname{Graph}_B(R;\vec{z})
	\end{equation}
	means that every polynomial equation $E(\vec{z})=0$ in the finite graph expansion is emitted as $sE(\vec{z})=0$.
\end{definition}

\begin{table}[ht]
	\centering
	\scriptsize
	\renewcommand{\arraystretch}{1.4}
	\begin{tabular}{@{}p{3.0cm}p{3.0cm}cp{3.4cm}@{}}
		\toprule
		\textsf{Constraint} & \textsf{Bounded expansion} & $\boldsymbol{\delta}$ & \textsf{Mechanism} \\
		\midrule
		Pair step
			& $\operatorname{Pair}_{B_k}(t_{i,r},h_{i,r},t_{i,r+1})$
			& 2 & Finite graph expansion \\
		Tail initialization
			& $t_{i,0} - s = 0$
			& 1 & Starts extraction chain \\
		Line readout
			& $f_i - h_{i,i} = 0$
			& 1 & Exposes fixed line via chain \\
		Boolean digits
			& $d_{i,\kappa} - d_{i,\kappa}^2 = 0$
			& 2 & Forces $d_{i,\kappa} \in \{0,1\}$ \\
		Non-adjacency
			& $d_{i,\kappa}\cdot d_{i,\kappa+1} = 0$
			& 2 & \textsc{Zeckendorf} constraint \\
		\textsc{Zeckendorf} sum
			& $f_i - \sum_{\kappa=2}^K F_\kappa d_{i,\kappa} = 0$
			& 1 & Linear combination \\
		Axiom component
			& $E_{i,\alpha,q}:=A_{i,\alpha,q}(\vec{x})$
			& 2 & $\alpha\in\mathcal{J}^{\mathrm{ax}}_i$ \\
		MP component
			& $E_{i,\alpha,q}:=M_{i,\alpha,q}(\vec{x})$
			& 2 & Bounded MP expansion \\
		Selector booleanity
			& $s_{i,\alpha}(1-s_{i,\alpha})=0$
			& 2 & One active justification \\
		Selector totality
			& $\sum_{\alpha\in\mathcal{J}_i}s_{i,\alpha}-1=0$
			& 1 & Exactly one justification \\
		Guarded obligation
			& $s_{i,\alpha}E_{i,\alpha,q}(\vec{x})=0$
			& 3 & Linear selector $\times$ quadratic \\
		Target check
			& $f_{n-1} - t = 0$
			& 1 & Linear equality \\
		\bottomrule
	\end{tabular}
	\caption{\label{tab:constraints}Degree bookkeeping for the pre-aggregation bounded-domain cubic constraint system $\mathcal{C}_{k,t}$. Every occurrence of $\mathrm{pair}$, $\mathrm{hd}$, $\mathrm{tl}$, or $\mathrm{line}$ is replaced by the finite bounded graph expansion of Definition~\ref{def:bounded-graph-expansion}. For each line $i$, the finite set $\mathcal{J}_i=\mathcal{J}^{\mathrm{ax}}_i\cup\mathcal{J}^{\mathrm{mp}}_i$ contains the axiom-schema alternatives and the bounded \emph{modus ponens} witness alternatives. A line justification is selected by boolean variables. For the active selector, each obligation $E_{i,\alpha,q}=0$ is emitted as $s_{i,\alpha}E_{i,\alpha,q}=0$.}
\end{table}

\begin{definition}[Bounded proof constraint system]
\label{def:proof-constraint-system}
	For fixed resource parameter $k$ and target $t\leq B_k$, define $\mathcal{C}_{k,t}$ to be the finite system of polynomial equations, as shown in Table \ref{tab:constraints}, expressing:
	\begin{enumerate}[label=\textnormal{(\roman*)}]
		\item line extraction from the carryless proof-list code, with every primitive-recursive symbol replaced by a finite bounded graph expansion;
		\item formula well-formedness and target equality;
		\item axiom-schema matching;
		\item modus ponens checking with bounded witness indices;
		\item selector constraints choosing exactly one active justification from the finite set $\mathcal{J}_i$ at each line.
	\end{enumerate}
	If $t>B_k$, the slice is declared unsatisfied, matching Definition~\ref{def:bounded-proof-predicate}.
\end{definition}

\begin{proposition}[Local degree bound]
\label{prop:global-degree-bound}
	Every polynomial in $\mathcal{C}_{k,t}$ has total degree at most $3$.
\end{proposition}

\begin{remark}
	The bounded graph expansions are finite systems of degree $\leq 2$. Line extraction, equality, booleanity, bounds, and structural decoding are degree $\leq 2$. Axiom and modus ponens obligations are expanded into degree-$\leq 2$ equations. A selector variable has degree $1$, so activating such an obligation gives degree at most $3$. No degree-$3$ obligation is multiplied by a selector.
\end{remark}

\begin{theorem}[Bounded checker--constraint equivalence]
\label{thm:checker-constraint-equivalence}
	For each fixed resource parameter $k$ and target $t\leq B_k$,
	\begin{equation}
		\operatorname{BoundedProof}_k(t)
		\quad\Longleftrightarrow\quad
		\exists\vec{x}\leq B_k\;\bigwedge_{P\in\mathcal{C}_{k,t}}P(\vec{x})=0.
	\end{equation}
\end{theorem}

\begin{remark}
	In the forward direction, a valid bounded proof supplies the proof-list code, all line-extraction witnesses, the finite graph selectors, the justification selectors, and the auxiliary \textsc{Zeckendorf} witnesses; these assignments satisfy every constraint in $\mathcal{C}_{k,t}$. In the reverse direction, any bounded satisfying assignment decodes to a finite proof trace whose lines satisfy the required axiom or modus ponens conditions and whose final line is the target.
\end{remark}

\section{Aggregation and Bounded Slices}
\label{sec:aggregation-and-slices}

\begin{lemma}[Two-channel decomposition]
\label{lem:two-channel}
	Every integer-coefficient polynomial $P$ can be written as
	\begin{equation}
		P=A-B,
	\end{equation}
	where $A$ and $B$ have nonnegative coefficients and
	\begin{equation}
		P(\vec{x})=0\iff A(\vec{x})=B(\vec{x}).
	\end{equation}
	Degree does not increase.
\end{lemma}

\begin{lemma}[Mixed-radix aggregation]
\label{lem:no-carry}
	Let $A_i,B_i$ be nonnegative channel values satisfying
	\begin{equation}
		0\leq A_i,B_i<B
	\end{equation}
	for a fixed base $B\geq 2$. Then
	\begin{equation}
		\sum_i B^i A_i=\sum_i B^i B_i
		\quad\Longleftrightarrow\quad
		\forall i\; A_i=B_i.
	\end{equation}
\end{lemma}

\begin{remark}
	This aggregation does not use sum of squares and therefore does not double the degree. The price is that digit bounds must be enforced. Here the enforcement is external: the theorem is a bounded-domain cubic satisfiability statement, with every variable ranging over $[0,B_k]$. The aggregation base is chosen larger than the maximum value of every channel $A_i,B_i$ on that finite domain. If one instead internalizes the bounds, the linear slack equations
	\begin{equation}
		A_i+u_i+1=B,\qquad B_i+v_i+1=B
	\end{equation}
	enforce $A_i,B_i<B$ over $\N$ and do not increase degree; this manuscript uses the external bounded-domain convention.
\end{remark}

\begin{definition}[Bounded cubic family]
\label{def:bounded-cubic-family}
	For each resource parameter $k$, the bounded construction emits a polynomial
	\begin{equation}
		U_k(u,x_1,\ldots,x_{N_k})\in\Z[u,x_1,\ldots,x_{N_k}]
	\end{equation}
	of total degree at most $3$. The variable count $N_k$, aggregation base, and representation bound are computable from $k$. The polynomial is interpreted over the bounded domain $[0,B_k]^{N_k+1}$.
\end{definition}

\begin{definition}[Bounded satisfaction]
\label{def:bounded-sat}
	For the slice indexed by $k$, define
	\begin{equation}
	\label{eq:sat-b}
		\mathrm{Sat}_{B_k}(t)
		:\iff
		\exists\rho\;\bigl(\rho(i_u)=t\;\wedge\;\mathrm{Bound}_{B_k}(\rho)\;\wedge\;U_k(\rho)=0\bigr),
	\end{equation}
	where
	\begin{equation}
		\mathrm{Bound}_{B_k}(\rho):\iff \forall i<N_k\;0\leq\rho(i)\leq B_k.
	\end{equation}
	Consequently $\mathrm{Sat}_{B_k}(t)$ is false when $t>B_k$, because the distinguished target coordinate is also bounded by $B_k$.
\end{definition}

\begin{proposition}[Decidability of bounded slices]
\label{prop:bounded-decidable}
	For each fixed $k$, the predicate $\mathrm{Sat}_{B_k}(t)$ is decidable.
\end{proposition}

\begin{proof}
	By exhaustive finite search over the bounded domain.
\end{proof}

\begin{theorem}[Aggregated bounded slice]
\label{thm:bounded-compilation}
	For each fixed resource parameter $k$ and every target code $t$,
	\begin{equation}
		\operatorname{BoundedProof}_k(t)
		\quad\Longleftrightarrow\quad
		\mathrm{Sat}_{B_k}(t).
	\end{equation}
	Both sides are false for $t>B_k$.
\end{theorem}

\begin{remark}
	This is the aggregated form of Theorem~\ref{thm:checker-constraint-equivalence}. The forward direction encodes a bounded proof as a bounded satisfying environment. The reverse direction decodes any bounded satisfying environment into a proof trace satisfying the bounded checker. Both directions are explicit.
\end{remark}

\begin{definition}[Slice-union satisfiability]
\label{def:slice-union-sat}
	Define
	\begin{equation}
	\label{eq:slice-union-sat}
		\mathrm{Sat}^{\cup}(t)
		:\iff
		\exists k\in\N\;\mathrm{Sat}_{B_k}(t).
	\end{equation}
	Equivalently,
	\begin{equation}
	\label{eq:slice-union-expanded}
		\mathrm{Sat}^{\cup}(t)
		\iff
		\exists k\in\N\;\exists\rho\;
		\bigl(\rho(i_u)=t\wedge\mathrm{Bound}_{B_k}(\rho)\wedge U_k(\rho)=0\bigr).
	\end{equation}
\end{definition}

\begin{remark}
	Definition~\ref{def:slice-union-sat} makes $\mathrm{Sat}^{\cup}$ a recursively enumerable predicate. It does \emph{not} assert the existence of a computable function $f:\N\to\N$ such that $\mathrm{Sat}^{\cup}(t)\iff \mathrm{Sat}_{B_{f(t)}}(t)$. The distinction is precisely the source of the error in earlier versions.
\end{remark}

\begin{proposition}[Computable bounding implies decidability]
\label{prop:computable-bound-decidable}
	If there were a computable function $f:\N\to\N$ such that
	\begin{equation}
		\mathrm{Sat}^{\cup}(t)\iff \mathrm{Sat}_{B_{f(t)}}(t),
	\end{equation}
	then $\mathrm{Sat}^{\cup}$ would be decidable.
\end{proposition}

\begin{proof}
	Given $t$, compute $f(t)$ and decide the fixed bounded slice $\mathrm{Sat}_{B_{f(t)}}(t)$ by finite search.
\end{proof}

\begin{remark}
	The standing result is a bounded cubic compilation theorem: for each fixed resource bound, a bounded-domain cubic slice is equivalent to bounded proof checking at that bound. The degree bookkeeping and the \textsc{Zeckendorf} encoding are unaffected by the corrigendum.

	The withdrawn claim required, in addition, a compression or uniformization principle that would turn the bounded family into a single fixed cubic Diophantine target for unbounded theoremhood. That gap is identified as the uniformity problem described in Section~\ref{sec:error}, and it remains open.
\end{remark}

\begin{theorem}[Corrected theorem]
\label{thm:corrected-theorem}
	For every fixed resource bound $k$ and every target code $t$ satisfying the explicit size condition $t\leq B_k$, there is an effectively computable finite system $\mathcal{C}_{k,t}$ of polynomial equations over $\mathbb{N}$, each of degree at most $3$, such that
	\begin{equation}
		\operatorname{BoundedProof}_k(t)
		\quad\Longleftrightarrow\quad
		\exists\vec{x}\leq B_k\;
		\bigwedge_{P\in\mathcal{C}_{k,t}}P(\vec{x})=0.
	\end{equation}
	Moreover, if the aggregation base is chosen larger than the maximum possible channel value on the bounded domain, this finite conjunction can be replaced by one cubic equation without increasing degree.
\end{theorem}

\begin{proof}
	Combine the local degree bound (Proposition~\ref{prop:global-degree-bound}), the bounded checker--constraint equivalence (Theorem~\ref{thm:checker-constraint-equivalence}), and the bounded aggregation theorem (Lemma~\ref{lem:no-carry}). Each claim has been established independently of the withdrawn uniformization step.
\end{proof}

\section{Cubic Thresholds}
\label{sec:hn-threshold}

\begin{remark}
	The corrected construction above is bounded and slice-wise. It should not be read as locating the classical degree-$3$ threshold for single Diophantine equations. Still, it is useful to record the elementary cubic consequence which follows from the known quartic constructions.
\end{remark}

\begin{lemma}[Quartic-to-cubic systems]
\label{lem:quartic-cubic-system-closure}
	Every quartic Diophantine equation in $v$ variables can be transformed, uniformly and effectively, into a finite system of Diophantine equations of degree at most $3$ in
	\begin{equation}
		v+\binom{v+1}{2}
	\end{equation}
	variables.
\end{lemma}

\begin{proof}[Proof sketch]
	Let the original variables be $x_1,\ldots,x_v$. For every pair
	$1\le i\le j\le v$, introduce a new variable $p_{ij}$ intended to represent the product $x_i x_j$, and add the quadratic equation
	\begin{equation}
		p_{ij}=x_i x_j.
	\end{equation}
	Every quartic monomial
	\begin{equation}
		x_a x_b x_c x_d
	\end{equation}
	can then be rewritten as a product of two quadratic product variables, for example
	\begin{equation}
		x_a x_b x_c x_d
		\quad\longmapsto\quad
		p_{ab}p_{cd},
	\end{equation}
	after fixing a uniform convention for ordering the pairs.

	The rewritten main equation has degree at most $2$ in the enlarged variables, while the defining equations $p_{ij}=x_i x_j$ have degree $2$. Cubic and lower degree terms of the original equation are left unchanged. Hence the output is a finite system of equations of degree at most $3$.

	Solvability is preserved. Any solution of the original quartic equation extends to a solution of the system by setting $p_{ij}=x_i x_j$. Conversely, any solution of the system satisfies the original quartic equation after substituting back the enforced products.
\end{proof}

\begin{corollary}[Trivial cubic family bound]
\label{cor:trivial-cubic-family-bound}
	Starting from \citeauthor{jones82}'s degree-$4$ universal equation in $58$ variables \citep{jones80,jones82}, one obtains a universal finite family of cubic equations using at most
	\begin{equation}
		58+\binom{59}{2}=1769
	\end{equation}
	variables.
\end{corollary}

\begin{proof}
	Apply Lemma~\ref{lem:quartic-cubic-system-closure} with $v=58$. The number of auxiliary product variables is
	\begin{equation}
		\binom{58+1}{2}
		=
		\binom{59}{2}
		=
		1711.
	\end{equation}
	Thus the total number of variables is
	\begin{equation}
		58+1711=1769.
	\end{equation}
\end{proof}

\begin{remark}
	The preceding subsection gives a coarse upper bound for cubic systems, but it does not locate the true single-equation boundary.
\end{remark}

\begin{definition}[Undecidability region]
\label{def:undecidability-region}
	Define
	\begin{equation}
		\mathcal{U} = \{(d,n)\mid \Hten{\N}{d}{n}=\mathsf{U}\}.
	\end{equation}
	For fixed $d$, define the threshold
	\begin{equation}
		n^*(d) = \min\{n\mid (d,n)\in\mathcal{U}\},
	\end{equation}
	when such an $n$ exists, and set $n^*(d)=\infty$ otherwise.
\end{definition}

\begin{proposition}[Arithmetical form of membership]
\label{prop:region-arithmetical-form}
	Fix a standard \textsc{Gödel Numbering} $e\mapsto\Phi_e$ of \textsc{Turing Machines} \citep{rogers67,kleene52}
	and a uniform coding of polynomial instances of type $(d,n)$. Let
	$\mathrm{Sol}_{d,n}(x)$ mean:
	\begin{quotation}
		\centering Instance $x$ has a solution in $\N$.
	\end{quotation}
	Then membership $(d,n)\in\mathcal{U}$ can be expressed arithmetically as the assertion that no machine is both total and correct for $\mathrm{Sol}_{d,n}$.
	Equivalently,
	\begin{equation}
		(d,n)\in\mathcal{U}
		\quad\Longleftrightarrow\quad
		\forall e\;
		\neg
		\Big(
			\mathrm{Tot}(e)
			\wedge
			\forall x\,
			\big(\Phi_e(x)=1\leftrightarrow \mathrm{Sol}_{d,n}(x)\big)
		\Big).
	\end{equation}
	In particular, the region $\mathcal{U}$ is not a primitive syntactic object; it is defined through a higher-level assertion about the nonexistence of deciders.
\end{proposition}

\begin{remark}
	Unpacking definitions: the problem $\Hten{\N}{d}{n}$ is undecidable exactly when no total \textsc{Turing Machine} computes its characteristic function \citep{turing37,turing38}. Totality is an arithmetical condition on a machine index, and correctness is expressed by comparing the machine output with the arithmetized solvability predicate $\mathrm{Sol}_{d,n}$. Thus membership in $\mathcal{U}$ is obtained by quantifying over all candidate deciders.
\end{remark}

\begin{proposition}[Monotonicity and threshold form]
\label{prop:region-monotonicity}
	The region $\mathcal{U}$ is upward closed in the parameters. If
	\begin{equation}
		(d,n)\in\mathcal{U},
	\end{equation}
	then for every $d'\ge d$ and $n'\ge n$,
	\begin{equation}
		(d',n')\in\mathcal{U}.
	\end{equation}
	Consequently, for each fixed degree $d$, the degree-slice of
	$\mathcal{U}$ has the form
	\begin{equation}
		\{n\mid n\ge n^*(d)\}
	\end{equation}
	if $n^*(d)<\infty$, and is empty otherwise.
\end{proposition}

\begin{proof}
	Enlarging the degree bound or the variable bound enlarges the input class. A decider for the larger class would restrict to a decider for the smaller class. Therefore undecidability propagates upward.
\end{proof}

\begin{remark}
\label{rem:boundary-obstruction}
	If the threshold function $d\mapsto n^*(d)$ were effectively available, then the undecidable region would have a computable boundary:
	\begin{equation}
		(d,n)\in\mathcal{U}
		\quad\Longleftrightarrow\quad
		n\ge n^*(d),
	\end{equation}
	with the convention that no $n$ qualifies when $n^*(d)=\infty$. This would not merely say that some Diophantine class is undecidable; it would determine the exact onset of undecidability for every degree. The known constructions do not provide such a boundary. They give isolated points, monotone consequences, and transformations between formats \citep{gasarch21,jones80}.

	This is the sense in which the region is impredicative: to know that a point lies in $\mathcal{U}$ is to know that every possible decision procedure fails on that entire class. The boundary is therefore described using the same computability-theoretic resources whose failure it records.
\end{remark}

\begin{remark}
\label{rem:region-relation-corrected}
	The bounded-domain cubic construction belongs on the safe side of this distinction. It gives explicit bounded slices and proves their local correctness. Each slice is decidable. The unbounded behavior appears only after taking the union over resource parameters:
	\begin{equation}
		\mathrm{Sat}^{\cup}(t)
		=
		\exists k\;\mathrm{Sat}_{B_k}(t).
	\end{equation}
	Thus the corrected result should be read as a bounded cubic compilation theorem.
\end{remark}

\begin{remark}
\label{rem:impredicativity-boundary}
	The degree-$3$ boundary remains delicate because several phenomena meet there: local proof checking, multiplication witnesses, selector activation, system aggregation, and the passage from bounded verification to unbounded existence. Each of these steps is individually natural, and each is correct in isolation. The error in the earlier versions was to treat their composition as automatically producing a single unbounded many-one reduction.

	Cubic systems are easily obtained from quartic equations. Bounded-domain cubic slices faithfully represent bounded proof checking. The open question is whether a compression principle exists that turns the bounded family into a single fixed cubic Diophantine target for unbounded theoremhood. We still hold that the question,
	\begin{quotation}
		At what point does the undecidability of undecidability begin?
	\end{quotation}
	is suggested to be, by the constructions shown, an instance of the very phenomenon it tries to pin down.
\end{remark}

\clearpage

\section*{Acknowledgments}
\vspace*{\fill}
{\scriptsize
	\bibliographystyle{plainnat}
	\setlength{\bibsep}{0.4em}
	\bibliography{refs}}
\vspace*{\fill}

\clearpage

\vspace*{\fill}

\subsection*{Revision History}

	This manuscript is the eighth and final release of our note. Revisions were driven by a precise boundary issue: bounded correctness, degree control, arithmetized proof checking, and unbounded completeness are closely related properties, but they are not interchangeable, and the passage between them is not automatic. The principal correction is the withdrawal of the \emph{many-one reduction} claim and its replacement with the slice-wise bounded compilation result, which was the valid core of the construction throughout.

\subsection*{Final Remarks}

	The author welcomes scholarly correspondence and constructive dialogue. No conflicts of interest are declared. This research received no funding.

\begin{center}
	\scriptsize{
	\vspace{3em}
	Milan Rosko
	\vspace{1em}
	ORCID: \href{https://orcid.org/0009-0003-1363-7158}{\textsf{0009-0003-1363-7158}}\\
	Email: \href{mailto:hi-at-milanrosko.com}{\textsf{hi-at-milanrosko.com}}\\
	\vspace{1em}
	Licensed under \enquote{Deed} \ccby \\ \href{http://creativecommons.org/licenses/by/4.0/}{\scriptsize\textsf{creativecommons.org/licenses/by/4.0}}
	}
\end{center}

\vspace*{\fill}

\end{document}